\documentclass{amsart}
\newif\ifcomments
\commentstrue  % change this to false to hide our rainbow comments
\usepackage{amsmath,amsthm,amssymb,amscd,url,enumerate}
\usepackage[utf8]{inputenc}
\usepackage{hyperref}
\usepackage{multirow}
\usepackage{pdflscape}
\usepackage{caption}
\usepackage[noadjust]{cite}
\usepackage{xypic}
\usepackage{ulem}
\usepackage{longtable}
\usepackage{tikz}
\usetikzlibrary{decorations.markings,arrows}
\usepackage{tikz-cd}
\usepackage{amsfonts}
\usepackage{amsmath}
\usepackage{graphicx}
\usepackage{amsthm}
\usepackage{xcolor}
\usepackage{mathrsfs}
\usepackage[margin=1in]{geometry}
\usepackage[ruled,vlined,linesnumbered,algosection]{algorithm2e}
\usepackage{algpseudocode}
\usepackage{csquotes}
\usepackage[all,cmtip]{xy} 
\usepackage{float}
\usepackage[caption = false]{subfig}
\usepackage{cleveref}

\usepackage{enumerate}
\usepackage{tcolorbox}
\theoremstyle{plain}
\newtheorem{theorem}{Theorem}[section]
\newtheorem{proposition}[theorem]{Proposition}
\newtheorem{lemma}[theorem]{Lemma}

\newtheorem{corollary}[theorem]{Corollary}
\newtheorem{example}[theorem]{Example}
\newtheorem{remark}[theorem]{Remark}
\newtheorem{definition}[theorem]{Definition}

\usepackage{listings}
\usepackage{color}

\newcommand{\iso}{\cong}

\DeclareMathOperator{\QQ}{\mathbb{Q}}

\SetEndCharOfAlgoLine{}
\SetKwInput{Input}{Input}
\SetKwInput{Output}{Output}
\SetKwFor{While}{While}{do}{}
\SetKwFor{ForEach}{For each}{do}{}
\SetKwIF{If}{ElseIf}{Else}{If}{then}{else if}{else}{endif}
\SetKw{Return}{Return}
\SetKw{KwRet}{return}
\SetKw{Recurse}{Recurse}
\SetKwFor{For}{For}{do}{}
\SetKw{Break}{break}

\newcommand{\QF}[1]{\mathbb{Q}(\sqrt{#1})}
\newcommand{\Z}{\mathbb{Z}}

\newcommand{\R}{\mathbb{R}}

\title{\textbf{Principal Well-Rounded Ideals of real quadratic fields }}

\author{Morgan Smith }
\address{Concordia University of Edmonton, 7128 Ada Blvd NW\\ Edmonton, Alberta,
T5B 4E3, Canada,}
\email{mamsmith0115@gmail.com}

\author{Ha T. N. Tran}
\address{University of Alberta -- Augustana Campus, 4901 -- 46 Avenue\\
Camrose, Alberta,T4V 2R3,  Canada,}
\email{htran2@ualberta.ca}

\keywords{Principal ideal, real quadratic field, Pell's equation, well-rounded ideal, lattice.}

\subjclass[2020]{11R11, 06B10, 11Y16, 11Y40, 11D09.}

\thanks{The authors acknowledge the support of the Natural Sciences and Engineering Research Council of Canada (NSERC) (funding RGPIN-2019-04209 and DGECR-2019-00428). We would like to thank the Center for Innovation and Applied Research (BMO-CIAR) for providing us with a computer for our calculations at the early stage of the project. We also would like to thank Gaurish Korpal for his useful comments on the first manuscript.
}

\begin{document}
	\maketitle
	
\begin{abstract}
It has been well known since Gauss that the principality of an ideal in a real quadratic field $K$ is equivalent to the solvability of a certain generalized Pell equations.  
In this paper, we combine this classical result with Srinivasan's conditions for the existence of well-rounded ideals in $K$ to obtain necessary and sufficient criteria for a real quadratic field to have principal well-rounded (PWR) ideals. Using these criteria, we prove that there are infinitely many real quadratic fields that have PWR ideals. Moreover, these ideals are pairwise non-similar. We then construct new algorithms that produce these PWR ideals, especially when the field discriminant is large. Our algorithms run in sub-exponential time theoretically; however, they are very fast in practice by employing some commonly used probabilistic algorithms for testing squarefreeness.  Finally, we briefly consider criteria for the existence of prime PWR ideals and show that there are infinitely many real quadratic fields that have prime PWR ideals.
\end{abstract}
	
%%%%%%%%%%%%%%%%%%%%%%%%%%%%%%%%%%%%%%%%%%%%%%%%%%%%%%%%%%%%%%%
\section{Introduction}

Well-rounded (WR) ideals are closely connected to many important mathematical problems, such as the shortest vector problem, the sphere packing problem, and the kissing number problem for ideals of number fields \cite{M13}. They also have valuable applications in coding theory \cite{VLL13, GBKTKH16, GTKH16, DGALH18, DKAGKH21}. Principal ideals can be represented by a single generator, offering a short representation of these ideals. This feature is especially preferable for applications in coding theory or in cryptography where the degree of number fields used is large. This is the initial motivation for us to investigate principal WR (PWR) ideals. Another motivation to work on these ideals is from studying WR twists \cite{DK19, LTT22} of the ring of integers of a number field. It has been shown that if there is a PWR ideal with a totally positive generator,  then the ring of integers can be twisted to that ideal.

The discussion of PWR ideals in real quadratic fields was started by  Fukshansky et al. in \cite{fukshansky2}. They showed that there is a finite number of complex quadratic fields containing PWR ideals and suggested further research into the number of real quadratic fields that contain PWR ideals \cite[Question 2]{fukshansky2}. In \cite{GTKH16}, Gnilke et al. proved that, there exist infinitely many real quadratic fields $\QQ(\sqrt{d})$ that have PWR ideals for squarefree positive integers $d$  such that $d\equiv 1,3 \mod 4$. Srinivasan  \cite{S20} completed the answer to Fukshansky et al.'s question by showing that there do not exist any PWR ideals in real quadratic fields $\QQ(\sqrt{d})$ when $d \equiv 2 \pmod 4$.

For real quadratic fields, it has been known since Gauss (article 243 and [Section V, \cite{gauss}]) that an ideal is principal if and only if a certain generalized Pell equation is solvable. In this paper, we combine this classical result with Srinivasan’s conditions for the existence of WR ideals \cite{S20}  to derive the necessary and sufficient criteria for a real quadratic field to have PWR ideals. This result is presented in \Cref{maintheorem} and proven in Propositions \ref{Eq3mod4}, and  \ref{Eq1mod4}.

We write $d= d_1 d_2$ for some odd, squarefree integers $d_1, d_2$ such that $\gcd(d_1, d_2)=1$. Note that Gnilke et al.  \cite{GTKH16} only considered two particular cases of our results (\Cref{maintheorem}) that is the case when $d_2=d_1+2$ or $d_2= d_1 +4$. 
We show that there are more general families of infinitely many real quadratic fields that have PWR ideals (Theorems \ref{thm:thinfinited3} and \ref{thm:thinfinited1}). We also prove that any two of these ideals which are not from the same field are not similar, thus there exist infinitely many non-similar PWR ideals from real quadratic fields (see \Cref{Thm:nonsimilar}).

One of the conditions required in \Cref{maintheorem} is the solvability of some generalized Pell's equation involving $d_1, d_2$, or, equivalently, the principality of some certain ideals of $\QQ(\sqrt{d})$. In case the discriminant of the field $\QQ(\sqrt{d})$ is large, it is known that two problems: solving generalized Pell's equations in \Cref{maintheorem} and testing the ideal principality in $\QQ(\sqrt{d})$  are classically hard (see \Cref{sec:strategy} for more discussion). Therefore, we propose a new strategy to search for real quadratic fields and their PWR ideals without solving these two hard problems (see \Cref{sec:strategy}). We then apply the strategy to construct new algorithms  (Algorithms \ref{alg:d1d23mod4},  \ref{alg:d1d21mod4}
and \ref{alg:2d1d21mod4}) to produce real quadratic fields and their PWR ideals of large norm. Theoretically, these algorithms run in sub-exponential time since they require testing squarefree and this is the most time-consuming step in those algorithms. However, they are still much faster than current algorithms for solving generalized Pell's equations or the principal ideal problem (see \Cref{sec:strategy} for more details). In practice, they run quite fast. For example,  using SageMath \cite{sagemath}  we could produce real quadratic fields of discriminant approximately $10^{240}$ and PWR ideals of norm approximately $10^{120}$ in several seconds just with a normal laptop (see \Cref{ex:larged} for an illustration) while other algorithms to solve generalized Pell's equations or principal ideal problems may not work for this large discriminant. Indeed, with the same laptop,  to find a generator of the same PWR ideal above, SageMath ran more than $12$ hours without result, and Pari/GP did not finish the task after $10$ hours running. Moreover, we show that the pairs $(d_1, d_2)$ of the form given in our Algorithms \ref{alg:d1d23mod4},  \ref{alg:d1d21mod4}
and \ref{alg:2d1d21mod4} are squarefree with probability at least $64\%$ which is the same as the probability of a random integer to be squarefree (see  \Cref{prop:probsqfree}).

Finally, we show some sufficient conditions for the existence of prime PWR ideals from real quadratic fields (Corollaries \ref{corprime} and \ref{cordoubleprime}) and prove that there are infinitely many non-similar such ideals (see \Cref{prop:inf_primePWR}). Such a field must have the form $\QQ(\sqrt{3})$ or $\QQ(\sqrt{d})$ with $d\equiv 1 \mod 4$ for some positive, squarefree integer $d$ and if exits, this prime ideal is unique up to similarity of the corresponding lattices (see \Cref{prop:cond_primePWR}).

In this paper, we consider only primitive integral ideals because any non-primitive integral WR ideal in a real quadratic field can be factored as the form $t I$ where $t\in \Z$ and $I$ is also a primitive integral WR ideal \cite{S20}.

%\HT{Add about software and cite them}
We use Pari/GP \cite{PARI2} and SageMath \cite{sagemath} for our experiments. The code can be found at \cite{Morgan_github}. 

The structure of the paper is as follows. In \Cref{sec:bacground}, we recall some basic knowledge related to WR ideals and their properties. The necessary and sufficient conditions for a real quadratic field to have PWR ideals are presented in \Cref{sec:Pell_equation}. In \Cref{sec:strategy}, we show our results of our experiment and discuss our strategy for finding PWR ideals. We prove that there are infinitely many real quadratic fields containing PWR ideals and these ideals are non-similar in \Cref{sec:infinite_family}. We then present our algorithms to produce PWR ideals in \Cref{sec:algorithms} and discuss prime PWR ideals in \Cref{sec:primePWR}.

%%%%%%%%%%%%%%%%%%%%%%%%%%%%%%%%%%%%%%%%%%%%%%%%%%%%%%%%%%%%%%%

%%%%%%%%%%%%%%%%%%%%%%%%%%%%%%%%%%%%%%%%%%%%%%%%%%%%%%%%%%%%%%
\section{Background}\label{sec:bacground}

In this section, we introduce some definitions and notation
used in the next sections.

Let $d$ be a positive squarefree integer. Then we denote by $K=\QF{d}$ the real quadratic field of discriminant $\Delta_K$, where \begin{equation*}
    \Delta_K =\begin{cases}
         d, \text{ if } d \equiv 1 \pmod 4 \text { and }\\
         4d, \text{ if } d \equiv 2,3  \pmod{4}.
    \end{cases} 
\end{equation*}
The ring of integers of $K$ is $O_K=\Z[\delta]$ where \begin{equation*}
    \delta = \begin{cases}
     \sqrt{d}, \text{ if } d \equiv 2,3 \pmod 4 \text{ and,}\\
     \frac{1 + \sqrt{d}}{2}, \text{ if } d \equiv 1 \pmod 4.
    
\end{cases} 
\end{equation*}

Any ideal $I$ of $O_K$ is also called an ideal of $K$.  Each ideal $I$ of $K$ can be generated by two elements $\alpha, \beta \in O_K$ over the ring $O_K$. In other words, $I= \{ a \alpha+ b \beta: a, b \in O_K\}$. In this case, we also write $I=(\alpha, \beta)$. If $I$ can be generated by a single element in $\alpha \in O_K$ over $O_K$, that is, $I= \{ a \alpha: a\in O_K\}$, then $I$ is called a principal ideal. 
In case $I$ is generated by $\alpha, \beta \in O_K$ over $\Z$, that is, $I= \{ a \alpha+ b \beta: a, b \in \Z\}$, we write $I = \langle \alpha, \beta \rangle_\Z$. 

\begin{definition}
    \label{def:lattices}
    A subset $L$ of  $\mathbb{R}^2$ is called a  \textit{(full-rank) lattice}  in $\mathbb{R}^2$ if there exist two linearly independent vectors $b_1, b_2 \in \mathbb{R}^2$, called a \textit{basis} of $L$, such that $L = b_1 \mathbb{Z} \oplus b_2 \mathbb{Z}$.
     We also write $L=\langle b_1, b_2 \rangle_\Z $.  
  
\end{definition}

We first recall the following result.
\begin{lemma}\label{UniqueItegralbasis}
    \cite{Buchmann1995}(Proposition 2.5)  
    Let $K=\QF{d}$. A subset $I$ of $O_K$ is an ideal if and only if there exist integers $a,b,m$ such that \[I = \left\langle ma, m\frac{b + \sqrt{\Delta_K}}{2} \right\rangle_\Z,\]
    $a>0,m > 0$, $4a|(\Delta_k - b^2)$, and $-a<b\leq a$, if $a>\sqrt{\Delta_K}$, or $\sqrt{\Delta_K} - 2a < b < \sqrt{\Delta_K}$, if $a < \sqrt{\Delta_K}$. This representation of $I$ is unique.
\end{lemma}
%\HT{Add a remark here }
\begin{remark}
    The \textit{norm} of the ideal in \Cref{UniqueItegralbasis} is $N(I)= ma$. In this paper, we only consider \textit{primitive} ideals, hence $m=1$ and $N(I)=a$.
\end{remark}
We will denote the embeddings of $K$ into $\R$ by  $\sigma_1, \sigma_2$ where
\begin{align*}
\sigma_1,\sigma_2: K & \hookrightarrow \R\\
    \sigma_1(x+y\sqrt{d}) &= x + y\sqrt{d},\\
    \sigma_2(x+y\sqrt{d}) &= x-y\sqrt{d}.
\end{align*}

Define the function \begin{align*}
    \Lambda: K &\rightarrow \R^2\\
    \alpha &\mapsto \langle \sigma_1(\alpha), \sigma_2(\alpha)\rangle.
\end{align*}
When we apply this map to ideals in $K$ we obtain lattices in $\R^2$ (see \cite{fukshansky1, fukshansky2} for more information). We will use $\Lambda(I)$ to denote the lattice in $\R^2$, which is the image of the ideal $I$ in $K$ under $\Lambda$.

\begin{definition}
    Let $L$ be a lattice in $\R^2$. 
    \begin{itemize}
        \item The \textit{minimum} of $L$ is \[ |L| := \min\{ \|x \|^2 : x \in L, x \not= (0,0)\},\] where $\|x\|$ denotes the usual Euclidean norm or the length of vector $x$ in $\R^2$.
        
        \item The set of \textit{minimal vectors} of $L$ is \[S(L) := \{x \in L : \|x \|^2 = |L|\}.\] 
        For any lattice $L$ in $\R^2$, $|S(L)| \in \{2,4,6\}$. 
        
        \item The lattice $L \subseteq \R^2$ is called \textit{well-rounded (WR)} if $S(L)$ spans $\R^2$. In other words, $L$ has two $\R$-linearly independent shortest vectors. 

         \item Let $L$ be a WR lattice. Then, it has a basis consisting of two minimal vectors. This basis is called a \textit{minimal basis}  of $L$. 

        \item An ideal $I$ of $K$ is called a \textit{well-rounded ideal} (WR ideal) if $\Lambda(I)$ is WR. In that case, if $\Lambda(\alpha), \Lambda(\beta)$ is a minimal basis of $\Lambda(I)$ for some $\alpha, \beta \in I$, then we also call $\alpha, \beta$ a \textit{minimal basis of $I$}.

    \end{itemize}
\end{definition}

Srinivasan \cite{S20}  proved an equivalent condition for the existence of WR ideals in a real quadratic field.

\begin{proposition}\label{S20}\cite{S20}(Theorem 1.1\footnote{The original paper states $\sqrt{\frac{\Delta_K}{3}} < a < \sqrt{3\Delta_K}$, however, the case $a=\sqrt{\frac{\Delta_K}{3}}$ does occur with the WR ideal $(2,1-\sqrt{3}) \subseteq \QF{3}$ and the case $a=\sqrt{3\Delta_K}$ occurs with the WR ideal $(6,3-\sqrt{3}) \subseteq \QF{3}$.})
    Let $K $ be a real quadratic field with discriminant $\Delta_K$. A primitive ideal $I$ in the ring of integers is WR if and only if $I = \langle a, \frac{a-\sqrt{\Delta_K}}{2}\rangle_\Z$ for some positive integer $a$ that satisfies $\sqrt{\frac{\Delta_K}{3}} \leq a \leq \sqrt{3\Delta_K}$. Moreover, $a$ is the norm of $I$ and  $a|\Delta_K$. In addition, there does not exist a real quadratic field with $d \equiv 2 \pmod{4}$ that contains a WR ideal. 
\end{proposition}

\begin{definition}
    Two lattices $L_1$ and $L_2$ in $\mathbb{R}^2$ are called \textit{similar}, denoted $L_1 \cong L_2$,  if there exists a positive real number $\gamma$ and a $2\times2$ real orthonormal matrix $U$ such that $L_2 = \gamma U L_1$.
\end{definition}

\begin{definition}\label{defAngleofLattice}
     Let $L$ be a WR lattice in $\R^2$. There is a minimal basis $\{u,v\}$ of $L$ such that the angle between $u$ and $v$ is in the interval $[\pi/3,\pi/2]$. This angle is an invariant of the lattice and is called the \textit{angle of $L$}. We will denote it by $\theta(L)$.
\end{definition}

\begin{lemma}\cite{fukshansky2}\label{AngleSimilar} Two WR lattices $L_1,L_2 $ in $\R^2$ are similar if and only if $\theta(L_1) = \theta(L_2)$.
\end{lemma}

%%%%%%%%%%%%%%%%%%%%%%%%%%%%%%%%%%%%%%%%%%%%%%%%%%
\section{PWR ideals of real quadratic fields and Pell's equations} \label{sec:Pell_equation}

In this section, we will discuss necessary and sufficient conditions on a pair of positive integers $(d_1, d_2)$ such that the real quadratic field $K=\QQ(\sqrt{d})$ has PWR ideals where $d=d_1 d_2$ is odd (see \Cref{maintheorem}). The condition for the ideal principality is equivalent to the solvability of a certain generalized Pell equation (see \eqref{eq:main_thmlike_Pell}).  This fact has been known since Gauss in terms of biquadratic forms in his article 243, also see \cite[Section V]{gauss} and Ribenboim \cite[6.11] {ribenboim} for more details. For completeness, we include a proof (without the use of binary quadratic forms) of this result in \Cref{Eq3mod4} and \Cref{Eq1mod4} for the specific ideals we are working with.

\begin{theorem}\label{maintheorem}
       % Let $d = d_1d_2$ be an odd, positive squarefree integer such 
       Let $d $ be an odd, positive squarefree integer and $K=\QF{d}$.
        Then $K$ has PWR ideals if and only if $d= d_1 d_2$, $d_1 < d_2 \leq 3d_1$ and one of the following generalized Pell's equations  is solvable
\begin{equation}\label{eq:main_thmlike_Pell}
  k^2d_2 - \ell^2d_1 =
    \begin{cases}
       \pm 2 & \text{if $d\equiv 3 \mod 4$}\\
      \pm 4 & \text{if $d\equiv 1 \mod 4$}.
    \end{cases}  
\end{equation}
    Moreover, any  PWR ideal $I$  of $K$ must have the form $I_1$ or $I_2$ as below.
$$I_1 = PP_1 \hdots P_r  \text{ and } I_2 = PQ_1\hdots Q_s \text{ if } d\equiv 3 \mod 4, \text{ and } $$
$$I_1 = P_1 \hdots P_r \text{ and }  I_2 = Q_1\hdots Q_s \text{ if } d\equiv 1 \mod 4,$$
here $d=d_1 d_2$, 
$d_1 = p_1 \hdots p_r$ and $d_2 = q_1\hdots q_s$ are the prime factorizations of $d_1$ and $d_2$, and        $P, P_i$ and $Q_j$ are the unique prime ideal above $2, p_i$ and $q_j$, respectively, in $K$. \\
Moreover, the two ideals $I_1$ and $I_2$ are similar.
\end{theorem}

\begin{proof}
By \cite[Remark 1]{S20}, any PWR ideal $I$ of $K$ must have the norm $ 2 d_1$ or $2d_2$ in case $d\equiv 3 \mod 4$ and $d_1$ or $d_2$ in case $d\equiv 1 \mod 4$ for some squarefree, comprime  integers $d_1, d_2$, with $d= d_1 d_2$.  Thus, $I$ has the form $I_1$ or $I_2$ given in the theorem, and then the ideal factorization of $I_i$ follows. The two ideals $I_1$ and $I_2$ are similar by \Cref{lem:I12similar}. The rest of the statement is obtained by the results of Propositions \ref{Eq3mod4} and \ref{Eq1mod4}.
\end{proof}

\begin{remark}
The condition for an ideal to be principal leads to a generalized Pell's equation, for example, see Equation \ref{eq:Pell_eq}. We then cancel out the factor $d_1$ or $d_2$ on both sides of such an equation, resulting in generalized Pell-like equations, as in Equation \ref{eq:main_thmlike_Pell}. In this paper, we also refer to the latter as generalized Pell equations.
\end{remark}

First, we consider the case when $d \equiv 3 \pmod 4$.

\begin{proposition}
 \label{Eq3mod4}
    Let $d = d_1d_2 \equiv 3 \pmod 4$ and with the notation in \Cref{maintheorem}.  Then $I_1$ and $I_2$ are PWR ideals if and only if $d_1 < d_2 \leq 3d_1$ and the generalized Pell's equation below has some integer solution $k, \ell$ 
    \[k^2d_2 - \ell^2d_1 = \pm 2.\]
     Moreover, $I_i$ has a minimal basis $  d_i+\delta, d_i -\delta$. 
\end{proposition}

\begin{proof}
Since $d$ is squarefree, we can assume that $d_1<d_2$.

First, by \Cref{S20}, the ideal $I_1$ is WR if and only if $I_2$ is WR if and only if $d_1 < d_2 < 3d_1$.

Second, we will prove that $I_i= \langle 2d_i, d_i +\delta\rangle_\Z$ for $i \in \{1,2\}$. Let $J_i=\langle 2d_i, d_i +\delta\rangle_\Z$. Then  $2d_i\in \Z_{>0}$ and $8d_i| (4d - 16d_i^2)$, because $d_1$ and $d_2$ are odd. We also have $\sqrt{\Delta_K} - 4d_i < 2d_i < \sqrt{\Delta_K}$ if $d_i < 2\delta$. Hence $J_i$ is an ideal of $O_K$ by \Cref{UniqueItegralbasis} and $N(J_i)=2d_i$. 
  Thus $J_i = I_i$ because $I_i$ is the unique ideal of norm $2d_i$. 
  
Now the ideal $I_1=\langle 2d_1, d_1 +\delta\rangle_\Z$ is principal if and only if there exists an element $\alpha = 2od_1 + k(d_1 + \delta) \in I_1$, $o,k\in\Z$, such that \[|N(\alpha)| = N(I_1) = 2d_1.\] Then \[\alpha = 2od_1 + k(d_1 - \delta) = d_1 \ell  - k\delta\]   where $\ell=2o + k$. Hence 
\begin{equation}\label{eq:Pell_eq}
    |N(\alpha)| = |d_1^2 \ell^2 - k^2d| = 2d_1.
\end{equation}
Hence $I_1$ is principal if and only if there is a solution to the generalized Pell's equation $k^2d_2 - \ell^2d_1 = \pm2$.%

We have that $I_1 I_2$ is the unique ideal of norm $4d$ because the discriminant of $\mathbb{Q}(\sqrt{d})$ is $4d$ and the unique ideal factorization property. Moreover, $N((2\delta)) = |N(2\delta)|= 4d$.  Thus $I_1I_2 = (2\delta)$. Therefore, $I_1I_2$ is principal. In other words,  $I_1$ is principal if and only if $I_2$ is principal. 

 Now we will show that $I_i$ has minimal basis $  \{d_i+\delta, d_i -\delta\}$. Let $\Lambda_i=\Lambda(I_i)$ for $i=1,2$. We have $\Lambda(d_i +  \delta), \Lambda(d_i -  \delta)$ is a basis of $ \Lambda_i$  since $I_i= \langle 2d_i, d_i +\delta\rangle_\Z$. We will first show that $\Lambda(d_i - \delta) $ is shortest in $\Lambda_i$. Assume  by contradiction that there exists a nonzero element $\alpha = 2ad_i + bd_i - b\delta \in I_i$ such that $\|\Lambda(\alpha)\|< \|\Lambda(d_i - \delta)\|$. It implies that 
 \begin{align}\label{eq:length}
        (2a+b)^2d_i^2 + b^2d < d_i^2 + d.
    \end{align}
If $b\ne 0$, then it follows from the inequality in \eqref{eq:length} that $2a+b= 0$.  Therefore $b^2\ge 4$, and  \eqref{eq:length} implies that $3d< d_i^2$, which contradicts the WR condition (see \Cref{S20}). Thus one must have $b=0$ and  $(4a^2-1)d_i^2< d \le 3d_i^2$ by \Cref{S20}. This leads to  $a=0$ and hence $\alpha=0$, a contradiction. 

Since $2$ elements $\Lambda(d_i +  \delta) \in \Lambda_i$ and $\Lambda(d_i -  \delta) \in \Lambda_i$ have the same length, they are both shortest in $\Lambda_i$. Thus these two elements form a minimal basis of $\Lambda_i$. Therefore,   $d_i +  \delta$ and $d_i -  \delta$ is form a minimal basis for $I_i$.
\end{proof}

\begin{remark}
  The only case where $d_2 = 3d_1$ is when $d=3$, and the field $\QQ(\sqrt{3})$ has PWR ideals $(2,1-\sqrt{3})$ and $(6,3-\sqrt{3})$.  
\end{remark}

Now we will prove similar conditions to \Cref{Eq3mod4} for PWR ideals when $d \equiv 1 \pmod 4$.

\begin{proposition}
 \label{Eq1mod4}
      Let $d = d_1d_2 \equiv 1 \pmod 4$ and with the notation in \Cref{maintheorem}. Then $I_1$ and $I_2$ are PWR ideals if and only if $d_1 < d_2 < 3d_1$ and the generalized Pell's equation below has some integer solution $k, \ell$ 
   \[k^2d_2 - \ell^2d_1 = \pm 4.\]
    Moreover,  $I_i$ has a minimal basis $(d_i+\sqrt{d})/2, (d_i -\sqrt{d})/2$. 

\end{proposition}

\begin{proof}
Assume that $d_1<d_2$. Similar to the proof of Proposition \ref{Eq3mod4}, we can show that $I_i = \langle d_i, \frac{d_i + \sqrt{d}}{2}\rangle_\Z $ for $i=1,2$. Moreover, since $|N(-2\delta + 1)| = d = N(I_1 I_2)$, we have that  $I_1 I_2= (-2\delta + 1)$   principal. Therefore, $I_1$ is principal if and only if $I_2$ is principal. 
The ideal $I_1$ is WR if and only if $d_1 < d_2 < 3 d_1$ (by \Cref{S20}), which is equivalent to that $I_2$ is WR.

Finally,  $I_1=\langle d_1, \frac{d_1 + \sqrt{d}}{2}\rangle_\Z $ is principal if and only if there exists an element $\alpha = ad_1 +b\frac{d_1+\sqrt{d}}{2} \in I_1$ for some $a,b \in \Z$ such that such that $N(I_1)= |N(\alpha)|$. The latter equation and since $N(\alpha) = \left( \frac{2ad_1 + bd_1}{2}\right)^2-\frac{b^2d}{4} = \pm d_1$ lead to the generalized Pell's equation  $k^2 d_2 - \ell^2 d_1 = \pm 4$ in the proposition. 

Note that if this equation has a solution, then  $k$ and $\ell$ are either both odd or both even, because $d_1, d_2$ are both odd by assumption.
Thus, 
$\beta = \frac{ \ell + k }{2} \frac{d_1+\sqrt{d}}{2} +\frac{ \ell - k }{2} \frac{d_1-\sqrt{d}}{2}$ is in $I_1$, since $I_1=\langle d_1, \frac{d_1 + \sqrt{d}}{2}\rangle_\Z=  \langle \frac{d_1 + \sqrt{d}}{2}, \frac{d_1 - \sqrt{d}}{2}\rangle_\Z$. The element $\beta \in I_1$  has norm $\pm d_1 = \pm N(I_1)$. Hence, $I_1$ is principal.

Let  $(d_1,d_2,k,\ell)$ satisfying the conditions of the proposition and let $\Lambda_i= \Lambda(I_i)$ for $i=1, 2$. Then $\Lambda\left( \frac{d_i + \sqrt{d}}{2}\right), \Lambda\left( \frac{d_i - \sqrt{d}}{2}\right)$ is a basis of $ \Lambda_i$ and $\left\|\Lambda\left( \frac{d_i + \sqrt{d}}{2}\right)\right\|^2 = \frac{d_i^2 + d}{2}$.    
Similarly to the proof of \Cref{Eq3mod4}, there are no non-zero vectors in this lattice $\Lambda_i$ with squared lengths shorter than $\frac{d_i^2+d}{2}$. Therefore the minimum of $\Lambda_i$ is $ \frac{d_i^2 + d}{2}$ and $\frac{d_i + \sqrt{d}}{2}, \frac{d_i - \sqrt{d}}{2}$ is a minimal basis of $I_i$.
\end{proof}

\begin{lemma}\label{lem:I12similar}
    The ideals $I_1$ and $I_2$ as defined in \Cref{maintheorem}  are similar.
\end{lemma}
\begin{proof}
    Define $d_1,d_2$ as in \Cref{Eq3mod4} such that $I_1$ and $I_2$ are WR. 
      From the minimal bases $d_i + \delta, d_i- \delta$ of $\Lambda_i$ (see Proposition  \ref{Eq3mod4})  we can see that \[\frac{\delta}{d_1}\Lambda_1 = \frac{\delta}{d_1}\begin{bmatrix}
        d_1 + \delta & d_1-\delta \\
        d_1 - \delta & d_1+\delta
    \end{bmatrix}\Z^2 = \begin{bmatrix}
        d_2 + \delta & d_2-\delta \\
        d_2 - \delta & d_2+\delta
    \end{bmatrix}\Z^2 = \Lambda_2.\]
    Therefore $\Lambda_1\iso \Lambda_2$.
\end{proof}

%%%%%%%%%%%%%%%%%%%%%%%%%%%%%%%%%%%%%%%%%%%%%%%%%%%%%%%%%%%%%%%
%%%%%%%%%%%%%%%%%%%%%%%%%%%%%%%%%%
\section{Our experiment results and strategy}\label{sec:strategy}
One of the goals of our work is to construct classical algorithms to produce many PWR ideals of real quadratic fields and their bases for employing them in applications, for example, in coding theory \cite{VLL13, GBKTKH16, GTKH16, DGALH18, DKAGKH21}.

Let $d_1, d_2$ be squarefree integers such that $\gcd(d_1,d_2)=1$ and  $d=d_1d_2$. Then we say that the pair $d_1,d_2$ represents PWR ideals if they satisfy the conditions in \Cref{maintheorem}. 
There are several methods to identify whether a given pair $d_1,d_2$ represents PWR ideals or not. We remark that the case when $k=\ell=1$ (i.e. $d_2=d_1+2$ or $d_2=d_1+4$) was already considered in \cite{GTKH16} and hence is not studied here. 

Applying Propositions \ref{Eq3mod4} and \ref{Eq1mod4}, one solution to this question would be to solve the given generalized Pell's equations in these propositions. However, no classical polynomial-time algorithm exists to solve these equations, and many generalized Pell's equations do not even have a solution. For example, the equation $x^2 - 4y^2 = 3$ has no integer solutions. Additionally, even if $k$ and $\ell$ exist, they may be quite large. Indeed,  when $d \equiv 3 \pmod 4$, from \cite{ConradPell} we can bound for $k$ and $\ell$ as below
\begin{align*} 
    k \leq \frac{\sqrt{u}+1/\sqrt{u}}{\sqrt{2d_2}}, \text{ and }
\ell \leq \frac{\sqrt{u} + 1/\sqrt{u}}{\sqrt{2d_1}}  
\end{align*} where $u$ is the fundamental unit of $K$. In other words,  $u=x+y\sqrt{d}$ where $\{x,y\}$ is the smallest solution to the Pell equation $x^2 - dy^2 = 1$ (see \cite{Lenstra02} for more discussions about solving such Pell's equations). When $d \equiv 1 \pmod 4$, one has 
\begin{align*}
    k \leq \frac{\sqrt{u}+1/\sqrt{u}}{\sqrt{d_2}}, \text{ and } 
\ell \leq \frac{\sqrt{u} + 1/\sqrt{u}}{\sqrt{d_1}}.   
\end{align*}  There has been lots of research into bounds on $u$. The best bounds currently known on $u$ are
$\sqrt{\Delta_K-4}+ \sqrt{\Delta_K}]/2 \le |u| < \exp(1/2 \sqrt{\Delta_K} (1/2 \log {\Delta_K} +1))$ by \cite{KATAYAMA1994385}. Hence $k$ and $\ell$ can be very large even for small values of $d_1,d_2$. For example, if we take $d_1=7$ and $d_2=13$. Then $d=91$ and $\Delta_K = 364$. Thus we have that $u< 1.36 \times 10^{23}$, $k \leq 721406311805$ and $\ell \leq 98673170373$.  This approach is therefore not very practical.

The second approach to tell if $d_1$ and $d_2$ (when $d_1d_2 \equiv 3 \pmod{4})$ represent PWR ideals is to check instead that the unique ideal of norm $2d_1$ (or norm $2d_2$) is principal. 
Determining whether a given ideal is principal or not is called the principal ideal problem (PIP). There are multiple algorithms to solve the PIP, for example, a method developed by Buchmann \cite{BUCHMANN1987}. However, no polynomial-time classical algorithm has been developed to solve this problem. Assuming the Generalized Riemann Hypothesis, the run time of the best known classical algorithm for this problem is sub-exponential time $L(\frac{1}{2}, b)$ here $L(a,b)=\exp(b n^a (\log)^{1-a})$ where $n$ is the input size \cite{B89,V00}. Hallgren \cite{Hallgren} has developed a polynomial-time quantum algorithm to solve Pell's equations and the PIP for real quadratic fields. Biasse and Song \cite{BS16} also provided efficient quantum algorithms for solving the PIP in arbitrary degree number fields.

Thus, for our classical approach, it is helpful to find alternative methods of computing pairs $(d_1,d_2)$ that do not require solving a generalized Pell's equation or the PIP.

\begin{figure}
    \centering
    \includegraphics[scale=0.8]{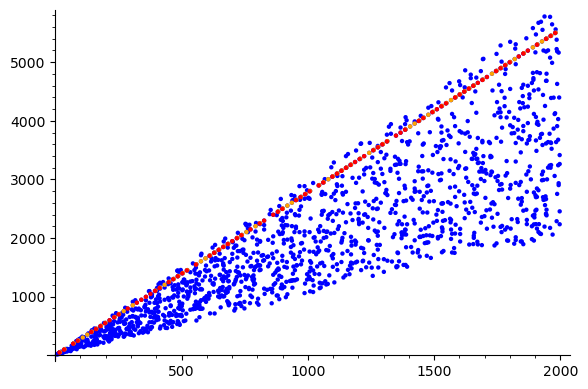}
    \caption{Values of $d_2$ found with $d_1<2000$ on the x-axis when $d \equiv 3 \pmod 4$ and $k < \ell < 8000$. The orange points are $d_2$ values which correspond to $k=3,\ell=5$ which $d_2k^2-d_1\ell^2 = 2$. The red points are $d_2$ values which correspond to $k=3,\ell=5$ with $d_2k^2-d_1\ell^2 = -2$.}
    \label{fig:3d1d2values}
\end{figure}

\begin{figure}\label{fi2}
    \centering
    \includegraphics[scale=0.8]{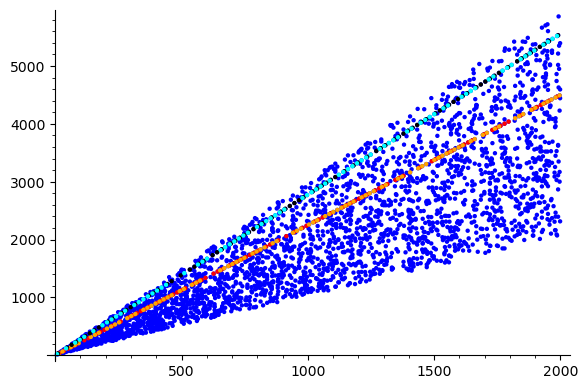}
    \caption{Values of $d_2$ found with $d_1<2000$ on the x-axis when $d\equiv 1\pmod4$ and $k < l \leq 50000$. The black points are $d_2$ values which correspond to $k=3,\ell=5$ with $d_2k^2-d_1\ell^2=4$. The cyan points are $d_2$ values which correspond to $k=3,\ell=5$ with $d_2k^2-d_1\ell^2 = -4$. The red points are $d_2$ values which correspond to $k=4,\ell=6$ with $d_2k^2-d_2\ell^2=4$. The orange points are $d_2$ values which correspond to $k=4,\ell=6$ with $d_2k^2-d_1\ell^2 = -4$.}
    \label{fig:1d1d2values}
\end{figure}

\vspace*{0.5cm}
\textbf{Numerical experiment:}
Using Pari/GP \cite{PARI2},  we generated  thousands of tuples $(d_1,d_2,k,\ell)$ which satisfy the generalized Pell's equations  in Propositions \ref{Eq3mod4} and \ref{Eq1mod4}. Some of these results are shown in Figures \ref{fig:3d1d2values} and \ref{fig:1d1d2values}. We performed an exhaustive search for $d_1<2000$ and $k,\ell <8000$ and with $d_1<10000$ and $k,\ell <50000$ when $d\equiv 3 \pmod{4}$ (see Figure \ref{fig:3d1d2values}), and for $d_1<10000$ and $k,\ell <50000$ in case $d\equiv 1 \pmod 4$ (see Figure \ref{fig:1d1d2values}). After investigating the obtained data, we found that there is a formula related to infinitely many such tuples $(d'_1, d'_2, k, \ell)$ to the smallest tuple $(d_1, d_2, k, \ell)$ as $d_1'=d_1 + 2k^2n, d_2'=d_2 + 2\ell^2n$ for some integer $n$ (see Theorems \ref{thm:thinfinited1} and \ref{thm:thinfinited3}).  This leads to our strategy as follows.

%\HT{add which axis  $d_1$ is, which one  $d_2$ is}

\vspace*{0.5cm}
%\newpage
\textbf{Our strategy:} We first consider the equations in Propositions  \ref{Eq3mod4} and \ref{Eq1mod4} as linear Diophantine equations instead of generalized Pell's equations. In other words, we choose some values of $k$ and $\ell$ first, then solve for one pair of $(d_1, d_2)$ satisfying Proposition \ref{Eq3mod4} or \ref{Eq1mod4}. Such a pair $(d_1, d_2)$ also gives us infinitely many other pairs $(d_1'=d_1 + 2k^2n, d_2'=d_2 + 2\ell^2n)$, for $n\in \Z$,  which also satisfy one of those propositions if they are squarefree.

%After having such a pair of $d_1$ and $d_2$,  one can obtain infinite many other pairs $(d_1', d_2')$ of the form $d_1'=d_1 + 2k^2n, d_2'=d_2 + 2\ell^2n$ for some integer $n$, which also satisfy Propositions \ref{Eq3mod4} and \ref{Eq1mod4} by checking whether   are squarefree or not. 

To illustrate our method, let's consider the case when $d= d_1 d_2 \equiv 3 \mod 4$. First,  we choose an odd integer $k$. Note that our method also works when $k$ is even given that $d \equiv 1 \pmod{4}$ (see Algorithm \ref{alg:2d1d21mod4}). Then we choose an integer $\ell$ such that $k < \ell < \sqrt{3}k$ and $\gcd(k, \ell)=1$. Next, using the extended Euclidean algorithm,  we solve for some positive integers $u$ and $v$ that satisfy the linear equation $k^2u -\ell^2 v = \pm 1$. Here $u$ and $v$ always exist since $\gcd(k, \ell)=1$. Now let $d_1= k^2 + 2v$ and $d_2= \ell^2 + 2u$. Then one can easily check that 
$k^2d_2- \ell^2 d_1=\pm 2$ is true (see the proof of Theorem \ref{thm:thinfinited3}). In other words, we found a pair $(d_1, d_2) $ that satisfies the generalized Pell's equation in Proposition \ref{Eq3mod4}. To obtain an initial pair of $d_1$ and $d_2$, we still need to test that both $d_1$ and $d_2$ are squarefree. After this initial pair, we can continue to generate more pairs $(d'_1, d'_2)$ by computing $d'_1= d_1 + 2k^2n$ and $d'_2= d_2 + 2\ell^2n$ and testing if they are squarefree for $n=1, 2, 3, \hdots$.

With the above strategy, we construct Algorithms \ref{alg:d1d23mod4}, \ref{alg:d1d21mod4}, and \ref{alg:2d1d21mod4} to produce PWR ideals of real quadratic fields, in particular when the discriminant of the field is large.   The most time-consuming step in these algorithms is to check whether $d_1$ and $d_2$ are squarefree. Theoretically, testing squarefree is still sub-exponential, it is however faster than solving generalized Pell's equations or solving the PIP. In particular, one can reduce testing squarefree to factoring which can be done in $L(\frac{1}{3}, b')$ by \cite{LL93} compared to $L(\frac{1}{2}, b)$ (for some constants $b, b'$) for solving  the PIP and  Pell's equations (hence generalized Pell's equations) by \cite{B89,V00}. In practice, testing squarefree is much faster using some probabilistic algorithms, see Example \ref{ex:larged} for more details. 

Finally, our algorithms can be easily adapted to quantum settings and then run in polynomial time, since factoring can be done in quantum polynomial time thanks to Shor \cite{Shor94}.

%%%%%%%%%%%%%%%%%%%%%%%%%%%%%%%%%%%%%%%%%%%%%%%%%%%%%%%%%%%%%%%

\section{The existence of infinitely many non-similar PWR ideals of real  quadratic fields} \label{sec:infinite_family}

In this section, we will prove that there are infinitely many real quadratic fields that have PWR ideals (see Theorems \ref{thm:thinfinited3} and   \ref{thm:thinfinited1}).  We do this by providing a series of lemmata showing that for given integers $k$ and $\ell$, satisfying certain criteria, there exist $d_1,d_2$ which satisfy the conditions of \Cref{maintheorem}  and that we can use this initial tuple $(d_1, d_2, k, \ell)$ to generate other such tuples. 
We also employ an invariant of WR ideals, their angle, to show that any two PWR ideals from two different real quadratic fields are non-similar. It follows that there are infinitely many non-similar PWR ideals as presented in \Cref{Thm:nonsimilar}.

First, we briefly recall the following lemma.
\begin{lemma}\label{Ricci} 
\cite{Ricci}
Let $f(x)$ be a separable polynomial function of degree 2 with integer coefficients. Suppose that $\gcd\{f(n):n\in \Z\}$ is a squarefree integer, then there are infinitely many squarefree values $f(n)$.
\end{lemma}

Now we consider the case that $d\equiv 3 \pmod{4}$.

\begin{lemma}\label{lemGen3mod4pt1}
    Let $k$ and $\ell$ be integers such that $k>0$, $k<\ell<3k$ and $\gcd(k,\ell)=1$. Then there exist integers $u,v$ such that $k^2u - \ell^2v = \pm 1$. Let $d_1 = k^2 + 2v + 2k^2n$ and $d_2 = \ell^2 + 2u + 2\ell^2n$ for some $n\in \Z$. Then $d_1 < d_2 < 3d_1.$
\end{lemma}

\begin{proof}
    Since the $\gcd$ of $k$ and $\ell$ is 1 we know that $u,v$ exist and can be found using the extended euclidean algorithm. 
     Consider the equation $k^2u-\ell^2v = \pm1$. 
    From this equation and $k<\ell$ we get \begin{align*}
        u &= \frac{\pm1+\ell^2v}{k^2} > \frac{\pm1 + k^2v}{k^2} = \frac{\pm1}{k^2} + v \geq v - 1, \text{and}\\
        3v &= \frac{3k^2u \mp 3}{\ell^2} > \frac{k^2u\mp1}{k^2} = u \mp \frac{1}{k^2} \geq  u - 1.
    \end{align*} Because the $\gcd$ of $u$ and $v$ must be $1$, one has $u>v$ and $3v > u$.
    Hence \begin{align*}
        d_2-d_1 &= \ell^2 + 2u + 2\ell^2n - k^2 - 2v - 2k^2n = (\ell^2 - k^2)(1+n)+2(u-v) > 0, \text{ and,}\\
        3d_1-d_2 &= 3k^2 + 6v + 6k^2n - \ell^2 - 2u - 2\ell^2n = (3k^2-\ell^2)(1+2n) + 2(3v-u) > 0.
    \end{align*}
    Thus $d_1<d_2<3d_1$.
\end{proof}

\begin{lemma}\label{lemGen3mod4pt2}
    Define $k, \ell, u, v, d_1, d_2,$ and $n$ as in Lemma \ref{lemGen3mod4pt1}, then $d =  d_1d_2 \equiv 3 \pmod 4$.
\end{lemma}

\begin{proof}
    We have that
    \begin{align*}
        d&=d_1d_2 = (k^2 + 2v + 2k^2n)(\ell^2 + 2u + 2\ell^2n) \\ 
        &\equiv (1+2v + 2n)(1+2u+2n) \pmod{4}. 
    \end{align*}
    Now, if $n$ is even $d \equiv (1 + 2v)(1+2u) \pmod{4}$ and if $n$ is odd $d \equiv (3+2v)(3+2u) \pmod{4}$. From the equation $k^2u-\ell^2v = \pm 1$,  we know that exactly one of $u,v$ is even because $k,\ell$ are both odd. Thus, \[d \equiv (1 + 2v)(1+2u) \equiv 1(1+2) \equiv 3 \pmod{4}\] when $n$ is even, and \[d \equiv (3 + 2v)(3+2u) \equiv 3(3+2) \equiv 3 \pmod{4}\] when $n$ is odd. Hence $d \equiv 3 \pmod{4}$ in both cases.
\end{proof}

\begin{lemma}\label{lemGen3mod4pt3}
    Define $k, \ell, u, v, d_1, d_2,$ and $n$ as in Lemma \ref{lemGen3mod4pt1}, then $d =  d_1d_2$ is squarefree for infinitely many values of $n$.
\end{lemma}

\begin{proof}
We will apply the result of  \Cref{Ricci} by assuming by contradiction that for some odd prime $p$, $p^2 | d$ for all $n \in \Z_{\geq 0}$. Then if we evaluate $d$ when $n=0,n=1$ and $n=2$ we have that $p^2$ divides $A, B, C$ where \begin{align*}
        A=& (k^2 + 2v)(\ell^2+2u) = k^2\ell^2 + 2v\ell^2 + 2uk^2 + 4uv,\\
        B=& (3k^2 + 2v)(3\ell^2+2u) = 9k^2\ell^2 + 6v\ell^2 + 6uk^2 + 4uv\text{ and }\\
        C=& (5k^2 + 2v)(5\ell^2+2u) = 25k^2\ell^2 + 10v\ell^2 + 10uk^2 + 4uv.
    \end{align*}
    
    Thus $p^2 | (C-2B +A) = 8k^2\ell^2$ and $p^2 | (3C-10B+15A) = 32uv$. Thus, $p^2| k^2 \ell^2$ and $p| uv$ since $p$ is odd. Hence, one has ($p^2|k^2$ or $p^2|\ell^2$) and ($p^2|u$ or $p^2|v$)  as $\gcd(k,\ell)=1$ and $\gcd(u,v)=1$.

    Now we have three cases to consider.
    
    \textbf{Case 1:} ($p^2 | k^2$ and $p^2 | v$) or ($p^2 | \ell^2$ and $p^2|u$).
    In this case, we have $p^2 | (k^2u - \ell^2v) = \pm1$, which is a contradiction.

    \textbf{Case 2:} ($p^2 | \ell^2$ and $p^2|v$).  
    If $p^2 | \ell^2$ and $p^2|v$, then 
    $$p^2| (k^2\ell^2 + 2v\ell^2 + 2uk^2 + 4uv-k^2\ell^2 -2k^2u - 4uv) = 2v\ell^2.$$ Thus $p | (k^2u - \ell^2v) =\pm 1$, which cannot happen.

    \textbf{Case 3:} ($p^2 | k^2$ and $p^2|u$). 
    Similarly to the previous case, if $p^2 | k^2$ and $p^2|u$, then 
    $$p^2| (k^2\ell^2 + 2v\ell^2 + 2uk^2 + 4uv-k^2\ell^2 -2\ell^2v - 4uv) = 2uk^2.$$  Thus $a^2 | uk^2$, which means that $p | \pm1$, a contradiction.
    
    Thus, no such $p$ exists and, by \Cref{Ricci}, there exist infinitely many squarefree values of $d$. 
\end{proof}

\begin{theorem}\label{thm:thinfinited3}
 For any two odd integers $k,\ell $ such that $\gcd(k,\ell) = 1$ and $k < \ell < \sqrt{3}k$,  there exist $d_1,d_2 \in \Z$, depending on $k$ and $\ell$,  which satisfy the conditions in Proposition \ref{Eq3mod4}.  

 Moreover, given any initial tuple $(d_1',d_2',k,\ell)$, which satisfies the conditions of Proposition \ref{Eq3mod4}, there exist infinitely many tuples $(d_1''=d_1' + 2k^2n, d_2''=d_2' + 2\ell^2n, k, \ell)$, $n \in \mathbb{Z}$, which also satisfy these conditions.% of Proposition \ref{Eq3mod4}. 
\end{theorem}

\begin{proof}
Suppose $k$ and $\ell$ are odd integers such that $k>0$, $k<\ell<\sqrt{3}\ell$ and $\gcd(k,\ell)=1$. Then there exist $g,h\in \Z$ such that $k^2g + \ell^2h = 1$, which can be found using the Extended Euclidean algorithm. We know $k,\ell > 0$, thus, either $g<0$ or $h<0$. Hence, we can write $k^2u - \ell^2v = \pm 1$ where $u=|g|$ and $v=|h|$.

    Let $n \in \Z_{\geq0}$, $d_1 = k^2 + 2v + 2k^2n$ and $d_2=\ell^2 + 2u + 2\ell^2n$. Then we have that \begin{align*}
        d_2k^2 - d_1\ell^2 &= (\ell^2 + 2u + 2\ell^2n)k^2 - (k^2 + 2v + 2k^2n)\ell^2 \\
        &= 2uk^2 - 2v\ell^2 = \pm 2.
    \end{align*}
     
    Thus, by Lemmas \ref{lemGen3mod4pt1}, \ref{lemGen3mod4pt2} and \ref{lemGen3mod4pt3}, we have shown that such $d_1, d_2$ satisfy all the conditions of \Cref{Eq3mod4} for some value of $n$. The first statement is then proven.

    Now we will prove the second statement of \Cref{thm:thinfinited3}. Suppose there exists a tuple $(d_1',d_2',k,\ell)$ that satisfies the conditions of  \Cref{Eq3mod4}.  Let $d''_1 = d_1' + 2k^2n$ and $d''_2=d_2' + 2\ell^2n$ for some positive integer $n$. We will prove that the tuple $(d_1'', d_2'', k, \ell)$ satisfies the conditions in \Cref{Eq3mod4}. First, it is easy to show that this tuple satisfies the general Pell equation in \Cref{Eq3mod4}. Second, we have that 
    \begin{align*}
        d''_1 d''_2= (d_1' + 2k^2n)(d_2' + 2\ell^2n)%\\ &= =d_1'd_2' + 2k^2d_2'n + 2\ell^2d_1'n + 4k^2\ell^2n^2 \\
        &\equiv 3 + 2n + 2n + 0 \equiv 3 \pmod 4.
    \end{align*}

    Third, we will show that $d_1''< d_2'' \leq 3d_1''$. Suppose $d_1'=1$. Then $d_2'=3$. Hence $\ell= -2 + 3k^2$ because $3k^2-\ell^2 = -2$ has no solution modulo $3$. Thus $d_2'' = d_2' + 2\ell^2n < 3(d_1' + 2k^2n) =3d_1''$. Alternatively, assume $d_1'>1$. Then \begin{align*}
        \ell^2 = \frac{\pm 2 + d_2'k^2}{d_1'} 
        < \frac{2 + 3d_1'k^2}{d_1'} < 3k^2 + 1
    \end{align*} 
    because $d_2'< 3d_1'$. Hence $\ell^2 \leq 3k^2$. Thus $d_2' + 2\ell^2n < 3(d_1' + 2k^2n)$. We also have that  $d_1' + 2k^2n < d_2' + 2\ell^2n$, for any values of $d_1'$ because $d_1' < d_2'$ and $k \leq \ell$. Therefore, $d_1'' < d_2'' < 3d_1''$.

    Now we will show that the following statement is true: for an infinite number of integers $n$, $d' = (d_1' + 2k^2n)(d_2' + 2\ell^2n)$ is squarefree. To prove this statement, we consider $d'$ as the function $f(x) = (d_1' + 2k^2x)(d_2' + 2\ell^2x)$. The two roots of $f(x)$ are $-\frac{d_1'}{2k^2}$ and  $-\frac{d_2'}{2\ell^2}$. These roots are distinct, otherwise, $d_1 '\ell^2 = d_2' k^2$, contradicts the fact that $d_2' k^2 - d_1' \ell^2 = \pm 2$.  
    Thus $f(x)$ is separable and of degree $2$. We also have that $f(0) = d_1'd_2'$ is squarefree by assumption. Thus $\gcd\{f(n) : n \in \mathbb{Z}\}$ is squarefree. Therefore, by Lemma \ref{Ricci}, the statement is held. 
   In other words, the tuple $(d_1'', d_2'', k, \ell)$ satisfies the conditions in \Cref{Eq3mod4} for infinitely many $n \in \Z_{>0}$.    
\end{proof}

\begin{remark}
    In \Cref{thm:thinfinited3}, $n$ can also be taken as a negative integer as long as $d_1d_2>0$.
\end{remark}

Now we will do similar to the above lemmata and \Cref{thm:thinfinited3} but for real quadratic fields $\QF{d}$ with $d\equiv 1\pmod{4}$.

\begin{lemma}\label{lemGen1mod4pt1}
    Let $k$ and $\ell$ be odd integers such that $k>0$, $k<\ell<3k$, and $\gcd(k,\ell)=1$. Then there exist integers $u,v$ such that $k^2u - \ell^2v = \pm 1$. Let $d_1 = k^2 + 4v + 2k^2n$ and $d_2 = \ell^2 + 4u + 2\ell^2n$ for some $n\in \Z$. Then $d_1 < d_2 < 3d_1$ and $d = d_1d_2 \equiv 1 \pmod{4}.$ Additionally, $d=d_1d_2$ is squarefree for infinitely many values of $n$.
\end{lemma}

\begin{proof}
    From \Cref{lemGen3mod4pt1} we know $u,v$ exist and $u>v$ and $3v > u$.
    Hence \begin{align*}
        d_2-d_1 &= \ell^2 + 4u + 2\ell^2n - k^2 - 4v - 2k^2n = (\ell^2 - k^2)(1+n)+4(u-v) > 0, \text{ and,}\\
        3d_1-d_2 &= 3k^2 + 12v + 6k^2n - \ell^2 - 4u - 2\ell^2n = (3k^2-\ell^2)(1+2n) + 4(3v-u) > 0.
    \end{align*}
    Thus $d_1<d_2<3d_1$.
     We also have that
    \begin{align*}
        d&=d_1d_2 = (k^2 + 4v + 2k^2n)(\ell^2 + 4u + 2\ell^2n) \\ 
        &\equiv (1+2n)(1+2n) \equiv 1 \pmod{4}. 
    \end{align*}
 The final statement of this lemma can be shown using a similar argument to the proof of \Cref{lemGen3mod4pt3}.
\end{proof}

\begin{lemma}\label{lemGen1mod4Evenpt1}
    Let $k$ and $\ell$ be even integers such that $k>0$, $\gcd(k,\ell)=2$, $8| (k\ell)$, and $k<\ell<3k$. Then there exist  integers $u,v$ such that $k^2u - \ell^2v = \pm 4$. If $u,v$ are both odd, let $n\in \mathbb{N}$, \begin{align*}
        d_1 &= \begin{cases}
            v + k^2(2n+1), \text{ if } u \equiv v \pmod{4}\\
            v + k^2(2n +1/2), \text{ if } u \not\equiv v \pmod{4},
        \end{cases}  \text{ and }\\
        d_2 &= \begin{cases}
            u + \ell^2(2n+1), \text{ if } u \equiv v \pmod{4}\\
            u + \ell^2(2n + 1/2), \text{ if } u \not\equiv v \pmod{4}.
        \end{cases}
    \end{align*}
    If $u$ or $v$ is even, let $u'$ be $u$ or $v$, respectively, and let $v'$ be $v$ or $u$, respectively. Let $n \in \mathbb{N}$,  \begin{align*}
        d_1 &= \begin{cases}
            v + k^2(n + 1/4), \text{ if } v' \equiv u' + 1 \pmod{4}\\
            v + k^2(n + 3/4), \text{ if } v' \equiv u' + 3 \pmod{4}, 
        \end{cases} \text{ and } \\
        d_2 &= \begin{cases}
            u + \ell^2(n + 1/4), \text{ if } v' \equiv u' + 1 \pmod{4}\\
            u + \ell^2(n + 3/4), \text{ if } v' \equiv u' + 3 \pmod{4}.
        \end{cases}
    \end{align*}
    Then $d_1 < d_2 < 3d_1$ and $d = d_1d_2 \equiv 1 \pmod{4}.$ Furthermore, $d$ is squarefree for infinitely many values of $n$.
\end{lemma}\begin{proof}
    One can prove this lemma by applying similar arguments to the proofs of \Cref{lemGen3mod4pt1} and \Cref{lemGen3mod4pt3}.

\end{proof}

\begin{theorem}\label{thm:thinfinited1}
For any two integers $k,\ell$ such that  $k < \ell < \sqrt{3}k$, and either
\begin{itemize}
        \item $k$ and $\ell$ are odd and  $\gcd(k,\ell)=1$, or 
        \item $k$ and $\ell$ are even,  $\gcd(k,\ell)=2$ and $8|k\ell$, 
\end{itemize}
there exist $d_1,d_2 \in \Z$, depending on $k, \ell$,  which satisfy the conditions of Proposition \ref{Eq1mod4}.

    Moreover, given an initial tuple $(d'_1,d'_2,k,\ell)$ that satisfies the conditions of \Cref{Eq1mod4},  there exist infinity many tuples of the form $(d_1'',d_2'',k,\ell)$, $n \in \mathbb{Z}$, which also satisfy these conditions where % of Proposition \ref{Eq1mod4}, where 
\begin{itemize}
    \item  $d_1'' = d_1' + 2k^2n $ and $d_2''= d_2' + 2\ell^2n$, if $k,\ell$ are odd, or 
    \item  $d_1''=d_1' + k^2n $ and $d_2''=d_2'+ \ell^2n$, if $k, \ell$ are even.
    \end{itemize}
    
\end{theorem} 

    \begin{proof}
        This can be shown using a similar argument to the proof of Theorem \ref{thm:thinfinited3} and by applying the results in Lemmas \ref{lemGen1mod4pt1} and \ref{lemGen1mod4Evenpt1}. 
    \end{proof}

\begin{theorem}   
\label{Thm:nonsimilar}
    There exist infinitely many real quadratic fields that have PWR ideals. Furthermore, any two of these ideals from distinct fields are non-similar. As a consequence, there are infinity many non-similar PWR ideals of real quadratic fields.
\end{theorem}
\begin{proof}
The statement that there are infinitely many quadratic fields with PWR ideals follows directly from Theorems \ref{thm:thinfinited3} and \ref{thm:thinfinited1}.

Define $d_1,d_2$ as in  \Cref{Eq3mod4} such that $I_1$ and $I_2$ are WR.     
  Since $\Lambda_1 \iso \Lambda_2$, we can consider the angle between the vectors of the minimal basis of $\Lambda_2$ which is $\{\Lambda( d_2 - \delta), \Lambda(d_2+\delta)\}$ by Proposition  \ref{Eq3mod4}. Then the cosine of this angle is 
  \[\frac{\Lambda(d_2 + \delta)  \cdot \Lambda( d_2 - \delta)}{\|\Lambda(d_2 + \delta) \| \cdot \|\Lambda(d_2 + \delta) \|} = \frac{d_2^2 - d}{d_2^2+d}=\frac{d_2-d_1}{d_2+d_1}.\] 
  Hence $0<  \frac{d_2-d_1}{d_2+d_1} \leq 1/2$. Therefore $\theta(\Lambda_2) = \arccos\frac{d_2-d_1}{d_2+d_1}=\theta(\Lambda_1)$ by Lemma \ref{AngleSimilar}.

    Let $d= d_1 d_2$, $c= c_1 c_2$ such that the pairs $(d_1,d_2)$, $(c_1, c_2)$ satisfy \Cref{maintheorem}. 
    Then the corresponding PWR ideals are similar if and only if \[\frac{d_2-d_1}{d_2+d_1} = \frac{c_2 - c_1}{c_2+c_1},\] by \Cref{AngleSimilar} and the above computation. Therefore $\frac{d_1}{d_2}= \frac{c_1}{c_2}$. 
    Then we must have $c_1=md_1$ and $c_2=md_2$ for some $m\in \mathbb{Q}$. Then there exit $ p,q \in \Z$ with $\gcd(p,q)=1$ such that $m=p/q$. Then $(p/q)d_1 \in \Z$ and $(p/q)d_2 \in \Z$. Thus $q|d_1$ and $q|d_2$, meaning $q=1$. Now $1=\gcd(c_1,c_2)=\gcd(pd_1,pd_2)=p$.  Hence $c_1= d_1,$ and $c_2= d_2$ and thus $d= d_1d_2= c_1 c_2= c$. Therefore, any PWR ideal lattices from different fields are not similar. %Thus there exist infinity many non-similar PWR ideals.
\end{proof}

It is known that among all two-dimensional lattices, the hexagonal lattice, denoted by $\mathcal{H}$, provides the highest density circle packing.  A result related to PWR ideals of real quadratic fields that are similar to the hexagonal lattice is as below.

\begin{corollary}\label{corHex}
    There exist exactly two primitive PWR ideals of real quadratic fields similar to the hexagonal lattice $\mathcal{H}$.  These ideals are $(2,1+\sqrt{3})$ and $(6,3+\sqrt{3})$ in $ \QF{3}$.
\end{corollary}
\begin{proof} 
Consider the ideal $I=(2,1+\sqrt{3})$ in $\QF{3}$. An integral basis of $I$ is $\{2, 1+\sqrt{3}\}$. Thus \[\Lambda(I) = \renewcommand{\arraystretch}{1}\begin{bmatrix}
        2&1+\sqrt{3}\\2&1-\sqrt{3} 
    \end{bmatrix}\Z^2 = \sqrt{2}\renewcommand{\arraystretch}{1}\begin{bmatrix}
        \frac{1}{\sqrt{2}}&\frac{1}{\sqrt{2}}\\\frac{1}{\sqrt{2}}&-\frac{1}{\sqrt{2}}
    \end{bmatrix}\mathcal{H}.\] 
    Hence $\Lambda(I)\iso \mathcal{H}$. From Proposition \ref{Eq3mod4} and Lemma \ref{AngleSimilar} we can also see that the ideal $(6,3+\sqrt{3})$ is PWR and $\Lambda((6,3+\sqrt{3}))$ is similar to $\mathcal{H}$. The uniqueness follows from \Cref{Thm:nonsimilar}.
\end{proof}

%%%%%%%%%%%%%%%%%%%%%%%%%%%%%%%%%%%%%%%%%%%%%%%%%%%%%%%%%%%%%%%

\section{ Algorithms to produce Principal Well-Rounded Ideals}
\label{sec:algorithms}

 In this section, by applying our strategy in \Cref{sec:strategy}, the results in \Cref{sec:infinite_family}, as well as the method of solving linear Diophantine equations,  we construct three algorithms, Algorithms \ref{alg:d1d23mod4}, \ref{alg:d1d21mod4} and \ref{alg:2d1d21mod4}, to produce PWR ideals of real quadratic fields. In addition, we will show that the probability a pair $(d_1, d_2)$ of the form in Step 6 of  Algorithms \ref{alg:d1d23mod4} and \ref{alg:d1d21mod4}, and in Step 22 of Algorithm \ref{alg:2d1d21mod4}, is squarefree is at least $64\%$ which is almost the same as the probability that a random integer is squarefree.

\begin{algorithm}
	\caption{Computing $d_1, d_2$ from $k$ for $d \equiv 3 \pmod 4$.  }\label{alg:d1d23mod4}
    \vspace{.2ex}
    \Input {An odd positive integer $k>1$.}

    \Output {Integers $d_1,d_2$ that satisfy the conditions in Proposition 
 \ref{Eq3mod4} or \ref{Eq1mod4}.}
     \vspace{.4ex}
     Choose an integer $\ell$ such that $k < \ell < \sqrt{3}k$ and $\gcd(k,\ell)=1$.\;
     Use the extended Euclidean algorithm to solve $k^2g + \ell^2h = 1$ for $g$ and $h$. \;
    
    Take $u=|g|$ and $v=|h|$ such that $k^2u-\ell^2v = \pm 1$.\;

    Set $d_1 \leftarrow k^2 + 2v$ and $d_2 \leftarrow \ell^2 + 2u$ and $n=0$. \;

    \While{$d_1$ or $d_2$ is not squarefree}{
    set $n \leftarrow n + 1$, $d_1 \leftarrow d_1 + 2k^2n$, $d_2 \leftarrow d_2 + 2\ell^2n$.   
     }

    \Return $d_1,d_2$. \;   
    \end{algorithm}

Here we will note that if we generate $d_1,d_2$ using Algorithm \ref{alg:d1d23mod4} we have that $n\ge 0$ and 
\begin{align*}
    \begin{cases}
        k^2 < d_1 < k^2(3/2 + 2n),\\
        k^2 < d_2 < 3k^2(3/2+2n)
    \end{cases}
\end{align*}

according to the bound of the coefficients generated using the extended Euclidean algorithm. 
One can also choose $n<0$ and then obtain 
 $d=d_1 d_2$  smaller than $k^4$. However, if $k$ is not sufficiently large,  it cannot be ensured that there exists a pair $(d_1,d_2)$ both are smaller than $k^4$ and squarefree.

We have performed Algorithm \ref{alg:d1d23mod4} with all $2<k<10000$ where $k<\ell<\sqrt{3}k$, $\gcd(k,\ell)=1$ and $d\equiv 3\pmod 4$. The largest value of $n$ required to find a squarefree pair $d_1, d_2$ was 9. In 70.77\% of cases the initial pair $d_1,d_2$ when $n=0$ was squarefree and in 21.35\% of cases we had $n=1$.  

The algorithm to find $d_1,d_2$ such that $d \equiv 1 \pmod 4$ differs from Algorithm \ref{alg:d1d23mod4} only in Step 4.
\begin{algorithm}
	\caption{\, Computing $d_1, d_2$ from $k$ for $d \equiv 1 \pmod 4$ with $k,\ell$ odd.  }\label{alg:d1d21mod4}
    \vspace{.2ex}
    \Input {An odd positive integer $k>1$.}

    \Output {Integers $d_1,d_2$ that satisfy the conditions in Proposition \ref{Eq3mod4} or \ref{Eq1mod4}.}
     \vspace{.4ex}
     Choose an integer odd $\ell$ such that $k < \ell < \sqrt{3}k$ and $\gcd(k,\ell)=1$.\;
     Using the Euclidean algorithm solve $k^2g + \ell^2h = 1$ for $g$ and $h$. \;
    
    Take $u=|g|$ and $v=|h|$ such that $k^2u-\ell^2v = \pm 1$.\;

    Set $d_1 \leftarrow k^2 + 4v$ and $d_2 \leftarrow \ell^2 + 4u$ and $n=0$. \;

    \While{$d_1$ or $d_2$ is not squarefree}{
    set $n \leftarrow n + 1$, $d_1 \leftarrow d_1 + 2k^2n$, $d_2 \leftarrow d_2 + 2\ell^2n$. 
    } 

    \Return $d_1,d_2$. \;   
    \end{algorithm}
    
    The bounds on $d_1,d_2$ for Algorithm \ref{alg:d1d21mod4} are similar to Algorithm 
 \ref{alg:d1d23mod4}. We have that \begin{equation*}
    \begin{cases}
        k^2 < d_1 < k^2(2+2n),\\
        k^2 < d_2 < 3k^2(2+2n).
    \end{cases}
\end{equation*}

\begin{algorithm}
	\caption{\, Computing $d_1, d_2$ from $k$ for $d \equiv 1 \pmod 4$ with $k,\ell$ even.  }\label{alg:2d1d21mod4}
    \vspace{.2ex}
    \Input {An even integer $k>2$}

    \Output {Integers $d_1,d_2$ that satisfy the conditions in Proposition 
 \ref{Eq1mod4}.} 
     \vspace{.4ex} 
     \If{$k$ is divisible by 4}{Choose an even integer $\ell$ such that $k < \ell < \sqrt{3}k$ and $4\not|\ell$.}
     \Else{Choose an even integer $\ell$ such that $k < \ell < \sqrt{3}k$ and $4|\ell$.}

Using the Euclidean algorithm solve $k^2g + \ell^2h = 4$ for $g$ and $h$. \;

Take $u=|g|$ and $v=|h|$ such that $k^2u-\ell^2v = \pm 4$.\;
    \If{$u,v$ are odd and $u \equiv v \pmod{4}$}{
Set $d_1 \leftarrow v$ and $d_2 \leftarrow u$ and $n=0$. \;

        }
    
    \If{$u,v$ are odd and $u \not\equiv v \pmod{4}$}{
Set $d_1 \leftarrow \frac{k^2}{2} + v$ and $d_2 \leftarrow \frac{\ell^2}{2} + u$ and $n=0$. \;

         }
    \If{$u$ is even}{
        \If{ $v \equiv u + 1\pmod{4}$}{
Set $d_1 \leftarrow \frac{k^2}{4} + v$ and $d_2 \leftarrow \frac{\ell^2}{4} + u$ and $n=0$. \;

            }
        \If{$v\equiv u + 3 \pmod{4}$}{
Set $d_1 \leftarrow \frac{3k^2}{4} + v$ and $d_2 \leftarrow \frac{3\ell^2}{4} + u$ and $n=0$. \;

            }}
    
    \If{$v$ is even}{
        \If{$u\equiv v + 1 \pmod{4}$}{
Set $d_1 \leftarrow \frac{k^2}{4} + v$ and $d_2 \leftarrow \frac{\ell^2}{4} + u$ and $n=0$. \;

             }
        \If{$u\equiv v + 3 \pmod{4}$}{
Set $d_1 \leftarrow \frac{3k^2}{4} + v$ and $d_2 \leftarrow \frac{3\ell^2}{4} + u$ and $n=0$. \;

            }}
            \While{$d_1$ or $d_2$ is not squarefree}{
    set $n \leftarrow n + 1$, $d_1 \leftarrow d_1 + k^2n$, $d_2 \leftarrow d_2 + \ell^2n$. 
    }
    \Return $d_1,d_2$. \; 
\end{algorithm}

Similarly to the case when $d\equiv 3 \pmod 4$, when we calculated $d_1,d_2$ for $d= d_1 d_2 \equiv 1 \pmod 4$ for $1 < k < 10000$ using Algorithm \ref{alg:d1d21mod4}, we found that 70.81\% of the results had $n=0,$ 21.60\% had $n=1$ and the largest $n$ value was 11.

Algorithm \ref{alg:2d1d21mod4} differs from the previous algorithms as it finds PWR ideals when the chosen integer $k$ is even. This requires that $d \equiv 1 \pmod{4}$.

\begin{remark}
     To find PWR ideals of real quadratic fields, it suffices to find a pair $d_1, d_2$ satisfying \Cref{maintheorem}. This can be done by applying Algorithms \ref{alg:d1d23mod4} and \ref{alg:d1d21mod4} and \ref{alg:2d1d21mod4}  which run in subexponential time in the worst case. These algorithms will be faster if we use some probabilistic algorithm for testing squarefree as in Pari/GP \cite{PARI2} or SageMath \cite{sagemath}. 
     
     By Theorems \ref{thm:thinfinited3} and \ref{thm:thinfinited1}, continuing to run the above algorithms will give us infinitely many values $d$, of which the field $\mathbb{Q}(\sqrt{d})$ has PWR ideals. In addition, for each algorithm, one can always choose $\ell=k+2$.
\end{remark}

\begin{example}\label{ex:larged}
    Applying the above algorithms and using SageMath \cite{sagemath}, we were able to compute\footnote{Calculations were run on a device with an Intel Core i7-1165G7 processor (2.80GHz, 4 cores) and 16 GB of RAM} PWR ideals when $d \approx 10^{240}$. To do this,  we first set $k=10^{60}-1$ and $l=k+2$. Running Algorithm \ref{alg:d1d23mod4},  we obtained the PWR ideals  $\langle d_i + \sqrt{d}, d_i-\sqrt{d} \rangle_\Z$ of the field $\mathbb{Q}(\sqrt{d})$ where $d= d_1 d_2 \equiv 3 \mod 4$ and  
    
    \begin{align*}    d_1=&19999999999999999999999999999999999999999999999999999999999955000 \\ &00000000000000000000000000000000000000000000000000000003, 
\end{align*} \begin{align*}    d_2=&200000000000000000000000000000000000000000000000000000000000350000 \\ &0000000000000000000000000000000000000000000000000000001\\
    &= d_1 + 7999999999999999999999999999999999999999999999999999999999998.
\end{align*} 
This calculation only took $4$ seconds.

Running Algorithm \ref{alg:d1d21mod4}, we obtained 
the PWR ideals  $\langle (d_i+\sqrt{d})/2, (d_i -\sqrt{d})/2\rangle_\Z$ of the field $\mathbb{Q}(\sqrt{d})$ where $d= d_1 d_2 \equiv 1 \mod 4$ and 
\begin{align*}    d_1=&499999999999999999999999999999999999999999999999999999999998900000 \\ &0000000000000000000000000000000000000000000000000000007,
\end{align*}  \begin{align*}
    d_2=&5000000000000000000000000000000000000000000000000000000000090000\\ &000000000000000000000000000000000000000000000000000000003\\
    &= d_1 + 19999999999999999999999999999999999999999999999999999999999996.
\end{align*}
This calculation took $632$ seconds (about $10.5$ minutes). It took longer than the calculation of the previous pair $d_1, d_2$ because the initial pair $d_1, d_2$ found was not squarefree. The code for computing this example can be found at  \cite{Morgan_github}.\\
Note that other algorithms to solve Pell's equations or principal ideal problems will not work for this large discriminant. We used the same laptop and ran the principal ideal test in SageMath \cite{sagemath} and it did not give any results even after $12$ hours. 
\end{example}

%%%%%%%%%%%%%%%%%%%%%%%%%%%%%%%%%%

We obtain the following result from \cite{Granville}(Theorem 1) in case $B=\gcd\{f(n): n \in \mathbb{N}\}$ is squarefree.

\begin{lemma}\label{lemGranville} Assume that the abc conjecture is true. Suppose that $f(x) \in \Z[x],$ is of degree two without any repeated roots. Let $B=\gcd\{f(n)| n\in \Z\}$ be squarefree. Then there are $\sim c_fN$ positive integers $n\leq N$ for which $f(n)$ is squarefree, where $c_f>0$ is a positive constant, which can be determined as follows:
\[c_f = \prod_{p \text{ prime}} \left( 1-\frac{w_f(p)}{p^2}\right),\]
where $w_f(p)$ is the number of integers $a$ in the range $1 \leq a \leq p^2$ for which $f(a) \equiv 0 \pmod{p^2}$. 
\end{lemma}

\begin{lemma}\label{lemwf1}
    Let $f(n)=d_1 (n) d_2(n)$ where $d_1(n)=k^2(1+2n)+gv$, $d_2(n)=\ell^2(1+2n)+gu$ with $g\in \{2,4\}$, and $k,\ell,u,v$ defined as in Algorithms \ref{alg:d1d23mod4} and \ref{alg:d1d21mod4}. Let $p$ be a prime. Then:
        \begin{itemize}
            \item $w_f(p)=0$ if $p=2$,
            \item $w_f(p)=1$, if $p|k$ or $p|\ell$, 
            \item otherwise $w_f(p)=2$. 
        \end{itemize}
        \end{lemma}
\begin{proof}
        For all values of $n\in \Z$, both $d_1(n)$ and $d_2(n)$ are odd integers. Hence $w_f(2)=0$.

        Let $p>2$ be a prime such that $p\nmid k$ and $p\nmid \ell$ and let $S=\{1,2,..,p^2\} \subseteq \Z$. Define the functions $h_1(n),h_2(n): S\rightarrow \Z_{p^2}$ such that $h_1(a)=d_1(a) \pmod{p^2}$ and $h_2(a)=d_2(a) \pmod{p^2}$. We will show that $h_1$ and $h_2$ are injective. Let $a,b \in S$ such that $h_1(a) \equiv h_1(b) \pmod{p^2}$. Then \[k^2(1+4a)+gv \equiv k^2(1+4b)+gv \pmod{p^2}.\] Thus $2ak^2 \equiv 2bk^2 \pmod{p^2}$. Hence $a \equiv b \pmod{p^2}$ because $p\nmid k$ and $p\nmid 2$. Then $a=b$ because $a,b \in S$. Therefore, $h_1$ is injective. Using the same argument, we can show that $h_2$ is also injective. Therefore, $h_1(a) \equiv 0 \pmod{p^2}$ for exactly one value $a \in S$ and $h_2(b) \equiv 0 \pmod{p^2}$ for exactly one value $b \in S$. Hence if $p\nmid k$, then $p^2|d_1$ for exactly one integer value of $n$ in the range $[1,p^2]$ and if $p\nmid \ell$, then $p^2|d_2$ for exactly one integer value of $n$ in the range $[1,p^2]$. Therefore, if $p\nmid k$ and $p \nmid \ell$, then $w_f(p)=2$ because $\gcd(d_1,d_2)=1$.

        Now we will consider when $p|k$ or $p|\ell$. Suppose that $p|k$. Then $p^2$ cannot divide $k^2(1+2n)+4v$ for any $n \in \Z$ because $\gcd(k,v)=1$. Similarly, if $p|\ell$, then $p^2$ cannot divide $\ell^2(1+2n)+4u$ for any $n \in \Z$ because $\gcd(\ell,u)=1$. Therefore, if $p|k$ or $p|\ell$, then $w_f(p)=1$ because $\gcd(k,\ell)=1$.
\end{proof}

\begin{lemma}\label{lemwf2}
    Define $k,\ell,u,v$ as in Algorithm \ref{alg:2d1d21mod4}. If $u,v$ are both odd, let \begin{align*}
        d_1(n) &= \begin{cases}
            v + k^2(2n+1), \text{ if } u \equiv v \pmod{4}\\
            v + k^2(2n +1/2), \text{ if } u \not\equiv v \pmod{4}, 
        \end{cases} \text{ and }\\
        d_2(n) &= \begin{cases}
            u + \ell^2(2n+1), \text{ if } u \equiv v \pmod{4}\\
            u + \ell^2(2n + 1/2), \text{ if } u \not\equiv v \pmod{4}.
        \end{cases}
    \end{align*}
    If $u$ or $v$ is even, let $u'$ equal $u$ or $v$, respectively, and let $v'$ equal $v$ or $u$, respectively. Let  \begin{align*}
        d_1(n) &= \begin{cases}
            v + k^2(n + 1/4), \text{ if } v' \equiv u' + 1 \pmod{4}\\
            v + k^2(n + 3/4), \text{ if } v' \equiv u' + 3 \pmod{4}, 
        \end{cases} \text{ and }\\
        d_2(n) &= \begin{cases}
            u + \ell^2(n + 1/4), \text{ if } v' \equiv u' + 1 \pmod{4}\\
            u + \ell^2(n + 3/4), \text{ if } v' \equiv u' + 3 \pmod{4}.
        \end{cases}
    \end{align*} Let $p$ be a prime and $f(n)= d_1(n) d_2(n)$. Then:
        \begin{itemize}
            \item $w_f(p)=0$ if $p=2$,
            \item $w_f(p)=1$, if $p|k$ or $p|\ell$, 
            \item otherwise $w_f(p)=2$.
        \end{itemize}
        \end{lemma}\begin{proof}
            This lemma can be proven using an argument similar to the proof of Lemma \ref{lemwf1}.
        \end{proof}

\begin{proposition}\label{prop:probsqfree}
  Fix two positive integers  $k$ and $\ell$. Let $f$ be defined as in Lemma \ref{lemwf1} or Lemma \ref{lemwf2}. Then, there are $\sim c_fN$ positive integers $n\leq N$ for which $f(n)$ is squarefree, where $c_f$ can be determined as follows:
   \[c_f = 2\prod_{p \text{ prime}}\left(1-\frac{2}{p^2} \right)\prod_{p>2 \text{ prime}, p|k\ell} \left(\frac{p^2-1}{p^2-2} \right).\]
   In particular $c_f > 0.64.$
\end{proposition}
\begin{proof}
    The form of $c_f$ follows from Lemmas \ref{lemGranville}, \ref{lemwf1} and \ref{lemwf2}.

    It is known that 
    $$\prod_{p \text{ prime}}\left(1-\frac{2}{p^2} \right) \approx 0.32263$$
    by \cite{OEIS}. Hence $2\prod_{p \text{ prime}}\left(1-\frac{2}{p^2} \right) \approx 0.64526$. Moreover, one has $\prod_{p>2 \text{ prime},  p|k\ell} \left(\frac{p^2-1}{p^2-2} \right) > 1.$
    Thus $c_f > 0.64.$
\end{proof}

%%%%%%%%%%%%%%%%%%%%%%%%%%%%%%%%%%%%%%%%%%%%%%%%%%%%%%%%%%%%%%%

\section{Prime Principal Well-Rounded  Ideals}\label{sec:primePWR}
In this section, let $K=\QQ(\sqrt{d})$ be a real quadratic field for some squarefree, positive integer $d$. 
We will briefly consider prime PWR ideals and sufficient conditions for their existence, after which we prove that there are infinitely many non-similar prime  PWR ideals. We show that such a field $K$ with $d>3$ has prime, PWR ideals only if $d \equiv 1 \mod 4$. If it exists, this prime ideal is unique up to similarity. %\HT{need to update}

\begin{proposition}\label{prop:cond_primePWR}
    Let $d>3$ be a positive, squarefree integer and $K= \QQ(\sqrt{d})$. If $K$ contains a prime WR ideal, then $d\equiv 1 \mod 4$, and the prime ideal is unique up to similarity.
\end{proposition}
        \begin{proof}
         From the proof of \Cref{Eq3mod4}, we know that any primitive PWR ideal in $K$ has norm $2d_1$ or $2d_2$ when $d\equiv 3 \mod 4$. Hence, no PWR ideal in $K$ is prime if $d \equiv 3 \mod 4$ and $d>3$.  Thus one must have that $d\equiv 1 \mod 4$. 
         
            For $K$ to contain a prime ideal we must have that $d=pm$ where $p$ is prime and $m\in \Z$ such that $p<m<3p$ or $m<p<3m$ by \Cref{S20}. Now we have three cases to consider.

           Case 1: $d=pq$ where $p$ and $q$ are distinct primes.
           
           If $d=pq$, $p<q<3p$, then clearly no non-similar WR ideals exist by  \Cref{S20}.
           
           Case 2: $d=p_1 \hdots p_n \cdot q$, $\prod_{i\le n}p_i < q < 3\prod_{i\leq n}p_i$ where all $p_i$ and $q$ are pairwise distinct. 
           Here we have that $p_i \geq 3$ for all $i\leq n$. Thus, $qp_j \geq 3q$, $0<j\leq n$. Hence $qp_j > 3(\prod_{i\leq n}p_i)/p_j$. The unique ideal of norm $qp_j$ (for $0<j\leq n)$ is not WR by \Cref{S20}.
           
           Case 3: $d=q \cdot p_1 \hdots p_n$, $q < \prod_{i<n}p_i < 3q$  where all $p_i$ and $q$ are pairwise distinct.           
           Pick $0< j \leq n$. Then $(\prod_{i \leq n}p_i)/p_j < q < p_jq$ because $p_j \geq 3$. Hence $p_jq>3(\prod_{i \leq n}p_i)/p_j$ because $q > (\prod_{i \leq n}p_i)/p_j$ and $p_j \geq 3$. The unique ideal of norm $p_jq$ is not WR by \Cref{S20}.
        \end{proof}
        
 \begin{remark}
     When $d=3$, we can take $d_1=1$ and obtain the prime PWR ideal $(2,1+\sqrt{3})$ of norm $2$. No other prime PWR ideals exist in this field.

 \end{remark}   

When $d \equiv 1 \pmod 4$ there exist prime PWR ideals in $\QF{d}$. For example, in $\QF{133}$ we have the prime PWR ideal $(7, \frac{7-\sqrt{133}}{2})$. These ideals occur when $d_1$ or $d_2$ is prime. One family of prime PWR ideals is easy to identify as below.

\begin{corollary}\label{corprime}
    Let $d=p(p\pm 4)$ where $p$ is a prime. If $p \pm 4$ is squarefree, then $K=\QF{d}$ has a prime PWR ideal.
\end{corollary}
    \begin{proof}
    This follows directly from  \Cref{Eq1mod4} with $k=\ell=1$.
\end{proof}

\begin{lemma}\label{walsh}\cite[Theorem]{trotter69} Let $p \equiv q \equiv 3 \pmod 4$, then $px^2-qy^2=\pm 4$ is solvable.
\end{lemma}

\begin{corollary}\label{cordoubleprime}
     Let $p$ and $q$ be primes. If $d=pq$, $p < q <3p$ and $p\equiv q \equiv 3 \pmod 4$, then $K=\mathbb{Q}(\sqrt{d})$ has a prime PWR ideal.
\end{corollary}
\begin{proof}
    This follows directly from \Cref{maintheorem},  \Cref{walsh} and \Cref{Eq1mod4}.
\end{proof}

   \begin{proposition}\label{prop:inf_primePWR}
       There exist infinitely many real quadratic fields containing prime PWR ideals and any two of these ideals from distinct fields are non-similar.
   \end{proposition}
    \begin{proof}
        From \cite[Theorem 2]{Mirsky}, it is known that for a large enough value of $x$, for any positive integer $H$, any non zero integer  $r$ and prime $p$, 
        \[|\{p\leq x: p-r \text{ is squarefree }\}| = \prod_{q \text{ is prime}, q\nmid r}\left( 1 - \frac{1}{q(q-1)} \right) \text{Li}(x) + o\left(\frac{x}{\log{x}^H}\right)\]
 where $\text{Li}(x)$ is the offset logarithmic integral. 
 Hence, there are infinitely many primes $p$ such that $p$ satisfies \Cref{corprime}. Therefore, there exist infinitely many real quadratic fields having prime PWR ideals.
The non-similarity of these ideals follows from \Cref{Thm:nonsimilar}.
    \end{proof}

%%%%%%%%%%%%%%%%%%%%%%%%%%%%%%%%%%%%%%%%%%%%%%%%%%%%%%%%%%%%%%%

%%%%%%%%%%%%%%%%%%%%%%%%%%%%%%%%%%%%%%%%%%%%%%%%%%%%%%%%%%%%%%%

%%%%%%%%%%%%%%%%%%%%%%%%%%%%%%%%%%%%%%%%%%%%%%%%%%%%%%%%%%%%%%%

\bibliographystyle{abbrv}
\bibliography{myrefs}

\begin{thebibliography}{10}

\bibitem{BS16}
J.-F. Biasse and F.~Song.
\newblock Efficient quantum algorithms for computing class groups and solving
  the principal ideal problem in arbitrary degree number fields.
\newblock pages 893--902, 2016.

\bibitem{B89}
J.~Buchmann.
\newblock A subexponential algorithm for the determination of class groups and
  regulators of algebraic number fields.
\newblock {\em S{\'e}minaire de th{\'e}orie des nombres, Paris},
  1989(1990):27--41, 1988.

\bibitem{Buchmann1995}
J.~Buchmann, C.~Thiel, and H.~Williams.
\newblock {\em Short Representation of Quadratic Integers}, pages 159--185.
\newblock Springer Netherlands, Dordrecht, 1995.

\bibitem{BUCHMANN1987}
J.~Buchmann and H.~Williams.
\newblock On principal ideal testing in algebraic number fields.
\newblock {\em Journal of Symbolic Computation}, 4(1):11--19, 1987.

\bibitem{ConradPell}
K.~Conrad.
\newblock Generalized {P}ell equation, {L}ecture {N}ote {M}ath 154, 2021.
\newblock
  \url{https://virtualmath1.stanford.edu/~conrad/154Page/handouts/genpell.pdf}.

\bibitem{DGALH18}
M.~T. Damir, O.~Gnilke, L.~Amor{\'o}s, and C.~Hollanti.
\newblock Analysis of some well-rounded lattices in wiretap channels.
\newblock In {\em 2018 IEEE 19th International Workshop on Signal Processing
  Advances in Wireless Communications (SPAWC)}, pages 1--5. IEEE, 2018.

\bibitem{DK19}
M.~T. Damir and D.~Karpuk.
\newblock Well-rounded twists of ideal lattices from real quadratic fields.
\newblock {\em Journal of Number Theory}, 196:168--196, 2019.

\bibitem{DKAGKH21}
M.~T. Damir, A.~Karrila, L.~Amorós, O.~W. Gnilke, D.~Karpuk, and C.~Hollanti.
\newblock Well-rounded lattices: Towards optimal coset codes for gaussian and
  fading wiretap channels.
\newblock {\em IEEE Transactions on Information Theory}, 67(6):3645--3663,
  2021.

\bibitem{fukshansky2}
L.~Fukshansky, G.~Henshaw, P.~Liao, M.~Prince, X.~Sun, and S.~Whitehead.
\newblock On well-rounded ideal lattices ii.
\newblock {\em International Journal of Number Theory}, 9(01):139--154, 2013.

\bibitem{fukshansky1}
L.~Fukshansky and K.~Petersen.
\newblock On well-rounded ideal lattices.
\newblock {\em International Journal of Number Theory}, 8(01):189--206, 2012.

\bibitem{gauss}
C.~F. Gauss.
\newblock {\em Disquisitiones arithmeticae}, volume 157.
\newblock Yale University Press, 1966.

\bibitem{GBKTKH16}
O.~W. Gnilke, A.~Barreal, A.~Karrila, H.~T.~N. Tran, D.~A. Karpuk, and
  C.~Hollanti.
\newblock Well-rounded lattices for coset coding in mimo wiretap channels.
\newblock In {\em 2016 26th International Telecommunication Networks and
  Applications Conference (ITNAC)}, pages 289--294. IEEE, 2016.

\bibitem{GTKH16}
O.~W. Gnilke, H.~T.~N. Tran, A.~Karrila, and C.~Hollanti.
\newblock Well-rounded lattices for reliability and security in rayleigh fading
  siso channels.
\newblock In {\em 2016 IEEE Information Theory Workshop (ITW)}, pages 359--363.
  IEEE, 2016.

\bibitem{Granville}
A.~Granville.
\newblock Abc allows us to count squarefrees.
\newblock {\em International Mathematics Research Notices}, 19:991--1009, 1998.

\bibitem{Hallgren}
S.~Hallgren.
\newblock Polynomial-time quantum algorithms for pell's equation and the
  principal ideal problem.
\newblock In {\em Proceedings of the Thiry-Fourth Annual ACM Symposium on
  Theory of Computing}, STOC '02, page 653–658, New York, NY, USA, 2002.
  Association for Computing Machinery.

\bibitem{KATAYAMA1994385}
S.~Katayama and S.~Katayama.
\newblock On bounds for fundamental units of real quadratic fields.
\newblock {\em Journal of Number Theory}, 46(3):385--390, 1994.

\bibitem{LTT22}
N.~H. Le, D.~T. Tran, and H.~T.~N. Tran.
\newblock Well-rounded twists of ideal lattices from imaginary quadratic
  fields.
\newblock {\em Journal of Algebra and Its Applications}, 21(07):2250133, 2022.

\bibitem{LL93}
A.~K. Lenstra and H.~W. Lenstra.
\newblock {\em The development of the number field sieve}, volume 1554.
\newblock Springer Science \& Business Media, 1993.

\bibitem{Lenstra02}
H.~W. Lenstra and F.~der Wiskunde~en Natuurwetenschappen.
\newblock Solving the pell equation.
\newblock pages 1--23, 2002.

\bibitem{M13}
J.~Martinet.
\newblock {\em Perfect lattices in Euclidean spaces}, volume 327.
\newblock Springer Science \& Business Media, 2013.

\bibitem{Mirsky}
L.~Mirsky.
\newblock The number of representations of an integer as the sum of a prime and
  a k-free integer.
\newblock {\em The American Mathematical Monthly}, 56(1):17--19, 1949.

\bibitem{OEIS}
{OEIS Foundation Inc.}
\newblock Decimal expansion of product $p$ prime $(1 - 2/p^2)$, entry {A}065474
  in the {O}n-{L}ine {E}ncyclopedia of {I}nteger {S}equences, 2025.
\newblock Published electronically at \url{http://oeis.org/A065474}.

\bibitem{ribenboim}
P.~Ribenboim.
\newblock {\em My numbers, my friends: Popular lectures on number theory}.
\newblock Springer Science \& Business Media, 2000.

\bibitem{Ricci}
G.~Ricci.
\newblock Ricerche aritmetiche sui polinomi.
\newblock {\em Rendiconti del Circolo Matematico di Palermo Series 1}, 57,
  1933.

\bibitem{Shor94}
P.~W. Shor.
\newblock Algorithms for quantum computation: discrete logarithms and
  factoring.
\newblock In {\em Proceedings 35th annual symposium on foundations of computer
  science}, pages 124--134. Ieee, 1994.

\bibitem{Morgan_github}
M.~Smith and H.~T.~N. Tran.
\newblock {PWR-Ideals-in-Real-Quadratic-Fields}.
\newblock
  \url{https://github.com/mamsmith/PWR-Ideals-in-Real-Quadratic-Fields}.

\bibitem{S20}
A.~Srinivasan.
\newblock A complete classification of well-rounded real quadratic ideal
  lattices.
\newblock {\em Journal of Number Theory}, 207:349--355, 2020.

\bibitem{PARI2}
{The PARI~Group}, Univ. Bordeaux.
\newblock {\em {PARI/GP version \texttt{2.13.4}}}, 2022.
\newblock available from \url{http://pari.math.u-bordeaux.fr/}.

\bibitem{sagemath}
{The Sage Developers}.
\newblock {\em {S}ageMath, the {S}age {M}athematics {S}oftware {S}ystem
  ({V}ersion 10.3)}, 2024.
\newblock \url{https://www.sagemath.org}.

\bibitem{trotter69}
H.~F. Trotter.
\newblock On the norms of units in quadratic fields.
\newblock {\em Proceedings of the American Mathematical Society},
  22(1):198--201, 1969.

\bibitem{VLL13}
R.~Vehkalahti, H.-F. Lu, and L.~Luzzi.
\newblock Inverse determinant sums and connections between fading channel
  information theory and algebra.
\newblock {\em IEEE transactions on information theory}, 59(9):6060--6082,
  2013.

\bibitem{V00}
U.~Vollmer.
\newblock Asymptotically fast discrete logarithms in quadratic number fields.
\newblock In {\em International Algorithmic Number Theory Symposium}, pages
  581--594. Springer, 2000.

\end{thebibliography}

%%%%%%%%%%%%%%%%%%%%%%%%%%%%%%%%%%%
\end{document}